\documentclass[12pt]{article}
\usepackage{amssymb}
\usepackage{latexsym}
\def\part#1{\frac{\partial\phantom{q}}{\partial#1}}

\newtheorem{thm}{Theorem}
\newtheorem{definition}{Definition}
\newtheorem{prp}[thm]{Proposition}
\newtheorem{lem}[thm]{Lemma}

\newcommand{\lie}[1]{\mathfrak{#1}}
\def\End{\mathop{\rm End}\nolimits}

\def\ker{\mathop{\rm ker}\nolimits}

\def\im{\mathop{\rm Im}\nolimits}

\def\tr{\mathop{\rm tr}\nolimits}

\def\i{\mathop{\iota}\nolimits}

\def\Diff{\mathop{\rm Diff}\nolimits}

\newcommand{\GCY}{generalized Calabi-Yau manifold\,\,}
\newcommand{\R}{\mathbf{R}}
\newcommand{\C}{\mathbf{C}}

\newcommand{\Z}{\mathbf{Z}}

\begin{document}
\title{Generalized Calabi-Yau manifolds}
 \author{Nigel Hitchin\\[5pt]
\itshape  Mathematical Institute\\
\itshape 24-29 St Giles\\
\itshape Oxford OX1 3LB\\
\itshape UK\\
 hitchin@maths.ox.ac.uk}
\maketitle
\begin{abstract}
\noindent A geometrical structure on even-dimensional manifolds is defined which generalizes the notion of a Calabi-Yau manifold and also a symplectic manifold. Such structures are of either  odd or even type and can be transformed by the action of both  diffeomorphisms and closed $2$-forms. In the special case of six dimensions we  characterize them as critical points of a natural variational problem on closed forms, and prove that a local moduli space is provided by an open set in either the odd or even cohomology.
 \end{abstract}

\section{Introduction} 
We introduce in this paper a geometrical structure on a manifold which generalizes both the concept of a Calabi-Yau manifold -- a complex manifold with trivial canonical bundle -- and that of a symplectic manifold. This is possibly   a useful setting for  the background geometry of   recent developments in string theory, but this was not  the original motivation for the author's first encounter with  this structure: it arose instead as part of a programme (following the papers \cite{H1},\cite {H2}) for characterizing special geometry in low dimensions  by means of invariant functionals of differential forms. In this respect, the dimension six is particularly important. This paper has two aims, then: first to introduce the general concept, and then to look at the variational and moduli space problem in the special case of  six dimensions.

We begin  with the definition in all dimensions of what we call  {\it generalized complex manifolds} and {\it generalized Calabi-Yau manifolds}\footnote{This terminology has been used before in \cite{Cand} but no confusion should arise in this paper.}. There are two novel features involved. The first is the use of the {\it Courant bracket}, a generalization of the Lie bracket on sections of the tangent bundle $T$ to sections of the bundle $T\oplus T^*$, and which comes to us from the study of constrained mechanical systems \cite{C}. The second is the {\it $B$-field} (this is the terminology of the physicists, but it is assuredly the same mathematical object). It turns out that the geometry we describe transforms naturally not only under the diffeomorphism group, but also by the action of a closed $2$-form $B$. 

To define a generalized complex manifold we imitate one definition of a K\"ahler   manifold. Instead of asking for the $(1,0)$ vectors to be defined by an isotropic  subbundle $E\subset T\otimes \C$, whose space of sections is closed under the  Lie bracket, we instead ask for a subbundle $E\subset (T\oplus T^*)\otimes \C$, isotropic with respect to the indefinite metric on $T\oplus T^*$ defined by the natural pairing between $T$ and $T^*$, and moreover whose space of sections is closed under the Courant bracket. For the definition of a generalized Calabi-Yau manifold we ask not for a closed $(n,0)$-form, but instead for a closed complex form $\varphi$ of mixed degree and of a certain algebraic type. This type is obtained by thinking of a form as a spinor for the orthogonal vector bundle $T\oplus T^*$ and then requiring the spinor to be {\it pure}. The well-known correspondence between maximally isotropic subspaces and pure spinors means that such a form defines a subbundle $E\subset (T\oplus T^*)\otimes \C$ and we show that sections of $E$ are closed under the Courant bracket if $d\varphi=0$.   There are two classes of such structures, depending on whether the degree of $\varphi$ is even or odd.

There are two motivating examples: an ordinary Calabi-Yau manifold and a symplectic manifold. A Calabi-Yau manifold with holomorphic $(n,0)$ form $\Omega$ defines a generalized Calabi-Yau structure by taking $\varphi=\Omega$. A symplectic manifold with symplectic form $\omega$ defines a generalized Calabi-Yau structure by taking $\varphi=\exp i\omega$. Transforming with a closed $2$-form $B$ means replacing $\varphi$ by $(\exp B)\wedge \varphi$. In certain cases, as we shall see, the B-field interpolates between symplectic and Calabi-Yau structures.

The special role of six dimensions arises from the fact that the group $\R^*\times Spin(6,6)$ has an open  orbit in either of its $32$-dimensional spin representations. Moreover a spinor in this open set  is the real part of a complex pure  spinor $\varphi$. We can define from this algebra  an invariant ``volume" functional defined on real forms, and we consider critical points of this functional on the closed                                                                                                               forms in a cohomology class in either the even or odd part of $H^*(M,\R)$ for a compact $6$-manifold $M$. If they lie in the open orbit  at each point of $M$, these critical points are precisely generalized Calabi-Yau manifolds. Imitating \cite{H1} we then show that, under a certain condition, a local moduli space for these structures is an open set in the corresponding cohomology group of even or odd degree. The required condition is implied by the $\partial\bar\partial$-lemma for complex manifolds and the strong Lefschetz theorem for symplectic ones. We should note that this approach forces us to consider two structures to be equivalent if they are related not just by the group of  diffeomorphisms isotopic to the identity, but by its extension by the action of   {\it exact} B-fields. 

There is a special pseudo-K\"ahler structure on the moduli space induced as a consequence of this approach. In the even case it is the structure determined by the intersection form -- ``without quantum corrections" in the physicists' language.

Finally, by returning to the origins of the Courant bracket, we observe that the whole structure can be twisted by a closed three-form, or  more naturally by a {\it gerbe} with connection.
\vskip .25cm
The author wishes to thank the Universidad Aut\'onoma, Madrid
 and the Programa Cat\`edra Fundaci\'on Banco de Bilbao y Vizcaya  for support during part of the preparation of this
paper.

\section{The Courant bracket}

We shall begin by setting up the less familiar pieces of differential geometry. The first is the bracket operation introduced by T. Courant (for $p=1$) in \cite{C}. This is an operation defined on pairs $(X,\xi)=X+\xi$ of a vector field $X$ and a $p$-form $\xi$ on a manifold $M$. Take  $X+\xi, Y+\eta \in C^{\infty}(T\oplus\Lambda^pT^*)$ and 
define
\begin{equation}
[X+\xi, Y+\eta]=
  [X,Y]+\mathcal{L}_{X}\eta -\mathcal{L}_{Y}\xi -\frac{1}{2}
  d(\i(X)\eta  -\i(Y)\xi)
  \label{courant}
  \end{equation}
  This operation is skew-symmetric but does not in general satisfy the Jacobi identity. One important feature  however is that, unlike the Lie bracket on sections of $T$, this bracket has non-trivial automorphisms defined by forms.  Let $\alpha \in\Omega^{p+1}$ be a closed $p+1$-form and 
  define the vector bundle automorphism $A$ of $T\oplus\Lambda^pT^*$ by 
  \begin{equation}
  A(X+\xi)=X+\xi+\i(X)\alpha
  \label{action}
  \end{equation}
  
 Then one may easily check that 
  $$A([X+\xi, Y+\eta])=[A(X+\xi),A(Y+\eta)].$$

\noindent{\bf Example:} Consider the case  $p=0$, so that $\xi$ is a function $f$. We then have 
$$[X+f, Y+g]=
  [X,Y]+Xg -Yf.$$
  This is the usual Lie bracket on $S^1$-invariant vector fields
  $$X+f\frac{\partial}{\partial\theta}$$
  on $M\times S^1$. Here the Jacobi identity does hold. Given a map $g:M\rightarrow S^1$, the diffeomorphism $(x,e^{i\theta})\mapsto(x,g(x)e^{i\theta})$ takes invariant vector fields to invariant vector fields and is the automorphism determined by the closed $1$-form $A=g^{-1}dg/i$.
\vskip .25cm
 We shall be concerned with the case $p=1$, where we have a bracket defined on sections of $T\oplus T^*$ and automorphisms defined by closed $2$-forms $B$. This interacts with a natural metric structure on $T\oplus T^*$ which we describe next.

  \section{The indefinite metric}
  
  \subsection{The orthogonal group}
  
  Let $V$ be an $n$-dimensional real vector space and consider the $2n$-dimensional space 
  $$V\oplus V^*.$$
  This admits a natural non-degenerate inner product of signature $(n,n)$ defined by 
$$(v+\xi,v+\xi)=-\langle v,\xi\rangle$$
for $v\in V$ and $\xi\in V^*$. Clearly the natural action of the general linear group $GL(V)$ preserves the inner product. The Lie algebra of the orthogonal group of \emph{all} transformations preserving the form  has the decomposition
 $${\mathfrak{so}}(V\oplus V^*)=\End V\oplus \Lambda^2V^*\oplus \Lambda^2V.$$
 In particular, $B\in \Lambda^2V^*$ acts as
 $$v+\xi\mapsto \i(v) B$$ and thus exponentiates to an orthogonal action on $V\oplus V^*$ given by
$$v+\xi\mapsto v+\xi+\i(v) B.$$
This is the algebraic action of a closed $2$-form which preserves the Courant bracket (\ref{action}).

\subsection{Spinors}\label{spinors}

Consider the exterior algebra $\Lambda^*V^*$ and the action of $v+\xi\in V\oplus V^*$ on it defined by
\begin{equation}
(v+\xi)\cdot\varphi=\i(v)\varphi+\xi\wedge \varphi
\label{Clifford}
\end{equation}
We have
$$(v+\xi)^2\cdot\varphi=\i(v)(\xi\wedge\varphi)+\xi\wedge\i(v)\varphi=(\i(v)\xi)\varphi=-(v+\xi,v+\xi)\varphi$$
which makes $\Lambda^*V^*$ into a  module over the Clifford algebra of $V\oplus V^*$. This defines the spin representation of the group $Spin(V\oplus V^*)$ if we tensor with the one-dimensional space $(\Lambda^n V)^{1/2}$. Splitting into even and odd forms we then have the two irreducible half-spin representations: 
\begin{eqnarray*}
 S^+&=&\Lambda^{ev}V^*\otimes (\Lambda^n V)^{1/2}\\ 
  S^-&=&\Lambda^{od}V^*\otimes
 (\Lambda^n V)^{1/2}
 \end{eqnarray*}
 If we now take $B\in \Lambda^2V^*\subset {\mathfrak{so}}(V\oplus V^*)$, then exponentiating $B$ to  $\exp B$ in the Lie group $Spin(V\oplus V^*)$ gives from (\ref{Clifford}) the following action on spinors:
$$\exp B (\varphi)=(1+B+\frac{1}{2}B\wedge B
+\dots)\wedge\varphi$$ 
\vskip .25cm
When $\dim V=n$ is even, there is an invariant bilinear form $\langle \varphi,\psi\rangle$ on $S^{\pm}$ which is symmetric if $n=4k$ and  skew-symmetric if $n=4k+2$. Using the exterior product we write this as 
$$\langle \varphi,\psi\rangle=\sum_m(-1)^m \varphi_{2m}\wedge\psi_{n-2m}\in \Lambda^{n}V^*\otimes ((\Lambda^n V)^{1/2})^2=\R$$
for even spinors and 
$$\langle \varphi,\psi\rangle=\sum_m(-1)^m \varphi_{2m+1}\wedge\psi_{n-2m-1}$$
for odd ones. If we define
$$\sigma:\Lambda^{ev/od}V^*\rightarrow \Lambda^{ev/od}V^*$$
by
$$\sigma(\varphi_{2m})=(-1)^m\varphi_{2m},\qquad \sigma(\varphi_{2m+1})=(-1)^m\varphi_{2m+1}$$
then
\begin{equation}
\langle \varphi,\psi\rangle=(\sigma\varphi\wedge\psi)_{n}.
\label{bilinear}
\end{equation}

\subsection{Pure spinors}

Given $\varphi\in S^{\pm}$, we consider its annihilator, the vector space 
$$E_{\varphi}=\{v+\xi\in V\oplus V^*: (v+\xi)\cdot\varphi=0\}$$
Since  $v+\xi\in E_{\varphi}$ satisfies
$$0=(v+\xi)\cdot(v+\xi)\cdot\varphi=-(v+\xi,v+\xi)\varphi$$
we see that $v+\xi$ is null and so  $E_{\varphi}$ is isotropic. A spinor $\varphi$ for which $E_{\varphi}$ is maximally isotropic (i.e. has dimension equal to $\dim V$) is called a {\it pure spinor}. Any two pure spinors are related by an action of $Spin(V\oplus V^*)$. 
To be pure is a non-linear condition which, in higher dimensions, is quite complicated. Here are some examples:
\vskip .25cm
\noindent{\bf Examples}:

\noindent 1. The spinor $1\in \Lambda^0 V^*\subset \Lambda^{ev}V^*$ is pure, since $(v+\xi)\cdot 1=\xi$ and so  the annihilator is defined by $\xi=0$, the maximal isotropic subspace $V\subset
V\oplus V^*$ 

\noindent 2. Applying any element of $Spin(V\oplus V^*)$ to $1$ gives another pure spinor. In particular we can exponentiate $B\in \Lambda^2V^*$ so that 
$$\exp{B}  =1+B+\frac{1}{2}B^2+\dots$$
is pure. Its  maximal isotropic
subspace is 
$$\exp B(V)=\{v+\i(v)B\in V\oplus V^*: v\in V\}$$
\vskip .25cm
The bilinear form applied to pure spinors has a geometrical meaning: $\langle\varphi,\psi\rangle=0$ if and only if $E_{\varphi}\cap E_{\psi}\ne 0$. For a proof of this see Chevalley \cite{Ch}, page 79.

Note that purity is scale-invariant, so we can define pure spinors not just in $S^{\pm}$ but also in the twisted spinor space $\Lambda^{ev/od}V^*$.

\section{Generalized complex structures}

\subsection{Definitions}
\begin{definition}
Let $M$ be a smooth manifold of dimension $2m$ with the indefinite metric on the bundle $T\oplus T^*$ defined by  $(v+\xi,v+\xi)=-\langle v,\xi\rangle$. A \emph {generalized complex structure} on $M$ is a subbundle $E\subset (T\oplus T^*)\otimes \C$ such that
\begin{itemize}
\item
 $E\oplus
\bar E=(T\oplus T^*)\otimes {\bf C}$
\item
the space of sections of $E$ is closed under the Courant bracket
\item
 $E$ is isotropic 
\end{itemize}
\end{definition}
The real version of this integrability -- a maximally isotropic subbundle of $T\oplus T^*$ with sections closed under Courant bracket -- is called a \emph{Dirac structure} in \cite{C}. A symplectic or Poisson structure on $M$ defines one of these. 
\vskip .25cm
Our main concern in this paper will be the notion  of a generalized Calabi-Yau manifold which we define next. Gualtieri's thesis \cite{MG} will contain more results on generalized complex manifolds.

\begin{definition}\label{GCYdef}
 A \emph {generalized Calabi-Yau structure} on a smooth manifold $M$ of dimension $2m$ is 
\begin{itemize}
\item 
 a closed form $\varphi\in \Omega^{ev}\otimes {\bf C}$ or $\Omega^{od}\otimes {\bf C}$ which is a complex pure
spinor  for the orthogonal vector bundle $T\oplus T^*$ and such that 
\item
  $\langle \varphi,\bar\varphi\rangle\ne 0$ at each point. 
  \end{itemize}
  
  \end{definition}
\vskip .5cm
The following proposition shows that a generalized Calabi-Yau manifold is a special case of a  generalized complex manifold.
\begin{prp} If $(M,\varphi)$ is a generalized Calabi-Yau manifold then the annihilator $E_{\varphi}\subset (T\oplus T^*)\otimes \C$ defines a generalized complex structure on $M$.
\end{prp}
\noindent{\bf Proof:} We saw from the algebra in the previous section that the annihilator of a pure spinor is maximally isotropic, so $E_{\varphi}$ certainly satisfies the last condition in the definition of generalized complex structure and has dimension $2m$. Moreover, since $\langle \varphi,\bar\varphi\rangle\ne 0$, we know that 
$$0=E_{\varphi}\cap E_{\bar\varphi}=E_{\varphi}\cap \bar E_{\varphi}$$
and so  
$$E_{\varphi}\oplus
\bar E_{\varphi}=(T\oplus T^*)\otimes {\bf C}.$$ 
It remains to show that sections of $E_{\varphi}$ are closed under the Courant bracket.

Suppose $X+\xi$ and $Y+\eta$ annihilate $\varphi$. Then from (\ref{Clifford})
$$\i(X)\varphi+\xi\wedge\varphi=0=\i(Y)\varphi+\eta\wedge\varphi$$
Using $d\varphi=0$  and ${\mathcal L}_X=d\i(X)+\i(X)d$ we obtain
\begin{eqnarray*}
\i([X,Y])\varphi&=&{\mathcal L}_X(\iota(Y)\varphi)-\i(Y){\mathcal L}_X\varphi\\
&=&-{\mathcal L}_X(\eta\wedge\varphi)-\iota(Y)d(\i(X)\varphi)\\
&=&-{\mathcal L}_X\eta\wedge \varphi-\eta\wedge {\mathcal L}_X\varphi+\i(Y)d(\xi\wedge\varphi)\\
&=&-{\mathcal L}_X\eta\wedge \varphi-\eta\wedge d(\i(X)\varphi)+\i(Y)(d\xi\wedge\varphi)\\
&=&-{\mathcal L}_X\eta\wedge \varphi+\eta\wedge d(\xi\wedge\varphi)+\iota(Y)(d\xi\wedge\varphi)\\
&=& -{\mathcal L}_X\eta\wedge \varphi+\eta\wedge d\xi\wedge\varphi+(\i(Y)d\xi)\wedge\varphi+d\xi\wedge\i(Y)\varphi\\
&=&-{\mathcal L}_X\eta\wedge \varphi+\eta\wedge d\xi\wedge\varphi+(\i(Y)d\xi)\wedge\varphi-d\xi\wedge\eta\wedge\varphi\\
&=&-{\mathcal L}_X\eta\wedge \varphi+(\i(Y)d\xi)\wedge\varphi
\end{eqnarray*}
and so, by skew symmetry,
\begin{eqnarray*}
\iota([X,Y])\varphi&=&\frac{1}{2}(\i([X,Y])\varphi-\iota([Y,X])\rho)\\
&=&\frac{1}{2}(-{\mathcal L}_X\eta\wedge \varphi+(\i(Y)d\xi)\wedge\varphi+{\mathcal L}_Y\xi\wedge \varphi-(\i(X)d\eta)\wedge\varphi)\\
&=&[\i(Y)d\xi+\frac{1}{2}d(\i(Y)\xi)-\i(X)d\eta-\frac{1}{2}d(\i(X)\eta)]\wedge\varphi\\
&=&[{\mathcal L}_Y\xi-{\mathcal L}_X\eta-\frac{1}{2}(d(\i(Y)\xi-\i(X)\eta)]\wedge\varphi
\end{eqnarray*}
From (\ref{courant}) this says that $[X+\xi,Y+\eta]\cdot\varphi=0$ as required.
\vskip .25cm
\noindent{\bf Examples:}

\noindent 1. Let $M$ be an $m$-dimensional complex manifold with a non-vanishing holomorphic form $\Omega$ of the top degree $m$. (A Calabi-Yau manifold is strictly speaking a K\"ahler manifold of this form -- there are others \cite{LT} -- but we shall use the terminology in the broader sense here.) Then $\Omega$ is pure since it is annihilated by any vector of the form $v+\xi$ where $v\in T\otimes \C$ is of type $(0,1)$ and $\xi\in T^*\otimes\C$ is of type $(1,0)$. This gives the maximal dimension $2m$ for the annihilator subspace. Since $\Omega$ has real  degree $2m$, the bilinear form is  $(-1)^m\Omega\wedge\bar\Omega\ne 0$. The algebraic conditions for a \GCY  are satisfied and since  $\Omega$ is closed  we obtain such a structure.

\noindent 2. Let $M$ be a symplectic manifold with symplectic form $\omega$. We saw that $1\in \Lambda^0T^*$ was pure, therefore, exponentiating the $2$-form $i\omega$, so is $\varphi=\exp i \omega  \in \Omega^{ev}\otimes \C$. The bilinear form gives
$$\langle\varphi,\bar\varphi\rangle=\frac{(-2i)^m}{m!}\omega^m$$
which is non-vanishing. Since $d\omega=0$, $\varphi=\exp i\omega$ defines a generalized Calabi-Yau manifold. 

\noindent 3. It is clear that the product of two generalized complex manifolds is a generalized complex manifold. Similarly if $(M_1,\varphi_1), (M_2,\varphi_2)$ are two generalized Calabi-Yau manifolds, then if $p_1,p_2$ denote the projections from the product $M_1\times M_2$,
$$\varphi=p_1^*\varphi_1\wedge p_2^*\varphi_2$$
defines a generalized Calabi-Yau structure on the product. The product of an odd type with an even type is odd and the product of two odd or two even types is even.

\subsection{The B-field}

If $B$ is a real closed $2$-form, and $(M,\varphi)$ a \GCY then 
$$(\exp B) \varphi=(1+B+\frac{1}{2}B^2+\dots)\wedge\varphi$$
is both closed and pure. Moreover, since $\exp B$ is real and acts through the Spin group 
$$\langle \exp B\,\varphi,\exp B \,\bar\varphi\rangle=\langle  \varphi,\bar\varphi\rangle\ne 0.$$
Thus a generalized Calabi-Yau structure  can be transformed by a B-field  to another one.

For a symplectic structure, the transform is 
$\varphi=\exp(B+i\omega)$
and we can always multiply $\varphi$ by a complex constant $c$, so one large class of generalized Calabi-Yau manifolds is given by
$$\varphi=c\exp(B+i\omega)$$
where $B$ is an arbitrary closed $2$-form and $\omega$ is a symplectic form.
\vskip .25cm
Suppose $M$ is a holomorphic symplectic manifold of complex dimension $2k$, for example a hyperk\"ahler manifold. It has a holomorphic non-degenerate $(2,0)$-form
$$\omega^c=\omega_1+i\omega_2.$$
The real and imaginary parts $\omega_1,\omega_2$ are themselves real symplectic forms on $M$. Let $t\ne 0$ be a real number then 
$$\exp (i\omega_2/t)$$
is the \GCY      defined by the symplectic form $\omega_2/t$. Apply the B-field  $B=\omega_1/t$ and we obtain
$$\exp ((\omega_1+i\omega_2)/t)$$
Multiply by the constant $t^k$ and we have a family of generalized Calabi-Yau structures defined by
$$\varphi_t=t^k\exp ((\omega_1+i\omega_2)/t)=t^k+\dots +\frac{1}{k!}(\omega_1+i\omega_2)^k.$$
Thus as $t\rightarrow 0$, these B-field transforms of the symplectic structure $\omega_2/t$ converge to the Calabi-Yau structure defined by the $(2k,0)$ form 
$(\omega_1+i\omega_2)^k/k!$.
In this way we can think of $B$ as interpolating between two extreme types of generalized Calabi-Yau manifold. 

\subsection{Dimension $2$}\label{two}

Let $M$ be a closed oriented surface. We consider the possible generalized Calabi-Yau structures on it.

First consider the {\it odd type}, defined by a form in $\Omega^1\otimes\C$. All such forms are pure. We thus have a closed complex $1$-form $\varphi$ such that
$$0\ne\langle\varphi,\bar\varphi\rangle=\varphi\wedge\bar\varphi.$$
This  non-vanishing  $1$-form is a $(1,0)$-form for a complex structure, and since it is closed, is holomorphic. This is therefore  an ordinary Calabi-Yau -- an \emph{ elliptic curve}. 
\vskip .25cm
 Now consider the {\it even type}. Here
$\varphi=c+\beta$, $c\in\Omega^0\otimes\C,\beta\in\Omega^2\otimes \C$. Since $\varphi$ is closed, $c$ is a  constant. Again, $\varphi$ is always pure, but we also have
\begin{equation}
0\ne \langle\varphi,\bar\varphi\rangle=c\bar\beta-\bar c\beta.
\label{ineq}
\end{equation}
In particular $c\ne 0$ and then from (\ref{ineq})
$$\beta/c-\bar\beta/\bar c=2i\omega\ne 0$$
so that $\omega$ is a symplectic form and 
$$\varphi=c\exp(B+i\omega)$$
where $2B=\beta/c+\bar\beta/\bar c.$
Thus the structure is the B-field transform of a symplectic manifold.

\subsection{Dimension $4$}\label{four}

On a $4$-manifold $M$ a generalized Calabi-Yau structure of odd type is defined by 
 $$\varphi=\beta+\gamma$$
 where $\beta$ is a complex closed $1$-form and $\gamma$ a complex closed $3$-form. The form $\varphi$ must define a complex pure spinor for $T\oplus T^*$. Here we are looking at the spin representation $S^-$ of the complexification $Spin(8,\C)$ of $Spin(4,4)$. In eight dimensions however, we have the special feature of triality -- the vector representation and the two spin representations are related by an outer automorphism of $Spin(8,\C)$. For us this means in particular that the  two spin spaces $S^{\pm}$ have the same structure as the vector representation -- an $8$-dimensional space with a non-degenerate quadratic form. The pure spinors are then just the null vectors in this space. 

It follows that $\varphi$ is pure if
\begin{equation}
0=\langle\varphi,\varphi\rangle=\beta\wedge\gamma.
\label{pure4}
\end{equation}
We also have the condition
\begin{equation}
0\ne \langle \varphi,\bar\varphi\rangle=\beta\wedge\bar\gamma+\bar\beta\wedge\gamma\ne 0
\label{nondegen4}
\end{equation}
which shows in particular that $\beta$ is nowhere vanishing. Thus from (\ref{pure4}), $\gamma=\beta\wedge\nu$ for some $2$-form $\nu$, well-defined modulo $\beta$. Using (\ref{nondegen4}) again,
\begin{equation}
\beta\wedge\bar\beta\wedge(\nu-\bar\nu)\ne 0
\label{nondegen42}
\end{equation}
and from this we can see that \emph{locally}, the structure on $M$ is defined by a map $f:M\rightarrow \C$ (where $df=\beta$) defining a fibration over an open set, a symplectic structure  $\Im \nu$ and a B-field $\Re \nu$ on the fibres. A global example is the product of an odd and an even $2$-dimensional generalized Calabi-Yau manifold.
 Tischler's theorem \cite{T} shows that a compact manifold with a non-vanishing closed $1$-form fibres over the circle and more generally that with two such forms like the real and imaginary parts of $\beta$, it must  fibre over $T^2$. In particular the first Betti number $b_1(M)$ is non-zero.
\vskip .25cm
For a structure of even type we have 
$$\varphi=c+\beta+\gamma$$
for a constant $c$, closed $2$-form $\beta$ and  $4$-form $\gamma$. For $\varphi$ to be pure we need
$$0=\langle\varphi,\varphi\rangle =2c\gamma-\beta^2.$$
If $c\ne 0$, this gives $\gamma=\beta^2/2c$.
The condition $0\ne\langle\varphi,\bar\varphi\rangle$ then gives
$$0\ne c\bar\gamma-\beta\bar\beta+\bar c \gamma=\frac{c\bar c}{2}(\bar\beta/\bar c-\beta/c)^2$$
so that $ \beta/c=B+i\omega$ where $\omega$ is symplectic. This gives
$$\varphi=c\exp(B+i\omega)$$
which is the transform of a symplectic structure.
\vskip .25cm
If $c=0$, the purity condition is $\beta^2=0$, which (as in the equation of the Klein quadric) means that $\beta$ is locally decomposable: $\beta=\theta_1\wedge\theta_2$. We also have
$$0\ne \langle\varphi,\bar\varphi\rangle=\beta\wedge\bar\beta=\theta_1\wedge\theta_2\wedge \bar\theta_1\wedge\bar\theta_2$$
so that $\theta_1,\theta_2$ span the space of $(1,0)$-forms for an almost complex structure and $\beta$ is of type $(2,0)$. Since $d\beta=0$ the structure is integrable and we have an ordinary Calabi-Yau manifold. In the compact case this must be a K3 surface or a torus. The remaining $4$-form $\gamma$ is the result of applying a (not necessarily closed) B-field to $\beta$.

\subsection{Structure groups and generalizations}

Our definition of a generalized complex structure yields a complex structure on $T\oplus T^*$ compatible with an indefinite metric. This is a reduction of the structure group of $T\oplus T^*$ to $U(m,m)\subset SO(2m,2m)$, together with an integrability condition. 

There are further reductions possible within this setting. First consider the case of a generalized Calabi-Yau manifold. Here the form $\varphi$ has the property $0\ne\langle\varphi,\bar\varphi\rangle$ so 
$$\frac{1}{\vert \langle\varphi,\bar\varphi\rangle\vert^{1/2}}\varphi$$
is a well-defined non-vanishing section of the spinor bundle, say
$$S^{+}\otimes \C\cong\Lambda^{ev}T^*\otimes (\Lambda^n T)^{1/2}\otimes \C.$$
Now $E=E_{\varphi}\subset (T\oplus T^*)\otimes \C$ is a complex maximally isotropic subbundle, so equally
$$S^{+}\otimes \C\cong\Lambda^{ev}E^*\otimes (\Lambda^n E)^{1/2}\otimes \C.$$
Here, the pure spinors with annihilator $E$ lie in $(\Lambda^n E)^{1/2}\otimes \C$, so $\varphi$ defines a non-vanishing section of $(\Lambda^n E)^{1/2}$ and, squaring it,  a trivialization of $\Lambda^n E^*$. This complex volume form on $E$  shows that a generalized Calabi-Yau structure  reduces the structure group to $SU(m,m)$.
\vskip .25cm
We could go further and consider smaller subgroups.  The four-dimensional case is instructive.
Here the spin representation defines a homomorphism (generating triality)
$$Spin(4,4)\rightarrow SO(4,4).$$
Restricting to the stabilizers of $2,3$ or $4$ spinors we obtain the special double covers: 
\begin{itemize}
\item
$SU(2,2)\rightarrow SO(4,2)$
\item
$Sp(1,1)\rightarrow SO(4,1)$
\item
$Sp(1)\times Sp(1)\rightarrow SO(4)$
\end{itemize}
The two spinors stabilized by $SU(2,2)$ are the real and imaginary part of a pure spinor  -- the complex closed form $\varphi$ of a generalized Calabi-Yau structure. The group $Sp(1,1)$ is obtained by requiring three forms to be closed giving a generalized \emph{hyperk\"ahler structure} and the last one involves four closed forms. The moduli space of such structures on a $K3$ surface has been studied by Nahm and Wendland \cite{NW}.

\section{The six-dimensional case}

\subsection{The quartic form}

We shall begin to study the algebra by working over the complex numbers, using $Spin(12,\C)$ instead of $Spin(6,6)$ and a complex six-dimensional vector space $V$. In this dimension the bilinear form on each of the $32$-dimensional spin spaces $S^{\pm}$ is skew symmetric, and so these are symplectic representations. The linear algebra we shall be doing is insensitive to the choice of orientation which distinguishes $S^+$ from $S^-$, but for various reasons the isomorphism
$$S^+\cong \Lambda^{ev}V^*\otimes (\Lambda^6 V)^{1/2}$$
will be a useful tool, so we shall fix $S=S^+$.
\vskip .25cm
 A symplectic action of a Lie group $G$ on a vector space $S$ defines a moment map 
 $$\mu:S\rightarrow {\mathfrak{g}}^*$$ given by
 $$\mu(\rho)(a)=\frac{1}{2}\langle\sigma(a)\rho,\rho\rangle$$
 where $\rho\in S$, $\sigma:{\lie g}\rightarrow \End S$ is the representation of Lie algebras and $a\in {\mathfrak{g}}$.
 On the Lie algebra $\mathfrak{so}(12,\C)$ we put the inner product $\tr XY$ and identify the Lie algebra with its dual, so $\mu(\rho)$ takes values in the Lie algebra. 
\vskip .25cm
\noindent{\bf Examples:}

\noindent 1.  Choose a basis vector $\nu$ for 
$(\Lambda^6 V)^{1/2}$ and consider the moment map for $Spin(12,\C)$ acting on $S$ at $\nu$. If $a=A+B+\beta$ in the decomposition  ${\mathfrak{so}}(V\oplus V^*)=\End V\oplus \Lambda^2V^*\oplus \Lambda^2V$, then 
$$\langle \sigma(a)\nu,\nu\rangle=\langle-(\tr A/2) \nu+B\wedge\nu,\nu\rangle=0$$
so the moment map vanishes on $\nu$ and hence on any pure spinor. 

\noindent 2. Now take $$\rho_0=\nu+\nu^{-1}\in (\Lambda^6V)^{1/2}\oplus (\Lambda^6V^*)(\Lambda^6V)^{1/2}\subset S.$$
In this case
$$\langle \sigma(a)\rho_0,\rho_0\rangle=-\tr A$$
and we find  that 
\begin{equation}
\mu(\rho_0)(v+\xi)=(-v+\xi)/4.
\label{momeq}
\end{equation}
\vskip .25cm
The moment map also defines an invariant:  
\begin{definition}
Let $\mu$ be the moment map for the spin  representation $S$ of $Spin(12,\C)$. Then
$$q(\rho)=\tr\mu(\rho)^2$$
is an invariant quartic function on $S$.
\end{definition}

This quartic has a close relationship with pure spinors:
\begin{prp}
\label{decomprop}
 For $\rho\in S$, $q(\rho)\ne 0$ if and only if $\rho=\alpha+\beta$ where $\alpha,\beta$ are pure spinors  and $\langle\alpha,\beta\rangle\ne 0$. The spinors $\alpha,\beta$ are unique up to ordering.
\end{prp}
\noindent{\bf Proof:} Consider as in the example  $\rho_0=\nu+\nu^{-1}\in S$: $\nu$ is pure with isotropic subspace $V$ and $\nu^{-1}$ with subspace $V^*$. 

Now suppose that $\alpha$ and $\beta$ are pure.  Because, up to a constant, $Spin(12,\C)$ acts transitively on pure spinors,  we can assume $\alpha=k\nu$. If  $\langle\alpha,\beta\rangle\ne 0$, we see from the definition of the bilinear form that $\beta_6\ne 0$ (we write $\alpha_{p}$ for the degree $p$ component of $\alpha$). By exponentiating an element of $\Lambda^2 V$ in the Lie algebra, we obtain a group element which leaves $\nu$ fixed but takes $\beta$ to an element 
$\tilde\beta$ with $\tilde \beta_4=0$ and 
 $\tilde \beta_6\ne 0$. But $\beta$ and hence $\tilde\beta$ are pure and so there is a $6$-dimensional isotropic space of vectors $v+\xi$ satisfying
$(v+\xi)\cdot\tilde\beta=0$. Looking at the degree $5$ term, this means that
$0=\i(v)\tilde\beta_6+\xi\wedge \tilde\beta_4=\i(v)\tilde\beta_6$ since $\tilde\beta_4=0$.
But then $v=0$ and the $6$-dimensional space is $V^*$. Thus $\tilde\beta=\ell\nu^{-1}$ and $\alpha+\beta$ can be transformed to 
$$k\nu+\ell\nu^{-1}.$$
From (\ref{momeq}) we see that  $q(\rho_0)=3$,  and so by invariance and homogeneity 
\begin{equation}
q(\alpha+\beta)= q(k\nu+\ell\nu^{-1})=3k^2\ell^2=3\langle\alpha,\beta\rangle^2
\label{qform}
\end{equation}
In particular,  the quartic invariant is non-zero for the sum of two pure spinors with $\langle\alpha,\beta\rangle\ne 0$.
\vskip .25cm
At $\rho_0=\nu +\nu^{-1}$ we saw that the moment map was $v+\xi\mapsto (-v+\xi)/4.$ Hence 
$\mu(\rho_0)^2=I/16$. If $\varphi=k\nu +\ell\nu^{-1}$ then $\mu(\rho)^2=k^2\ell^2I/16$ and so 
\begin{equation}
\mu(\rho)^2=\frac{1}{48}q(\rho)I
\label{square}
\end{equation}
By invariance this holds for all spinors in a $Spin(12,\C)$ orbit of $\rho_0$. 

Let $G\subset Spin(12,\C)$ be the stabilizer of $\rho_0$. Each element of $G$ 
 commutes  with $\mu(\rho_0)$ and so preserves or interchanges its 
 two eigenspaces $V,V^*$ in its action on $V\oplus V^*$. The identity 
 component $G_0$ preserves them, and hence lies in $GL(V)$. But if the group preserves $\nu^{-1}\in (\Lambda^6V^*)^{1/2}$ it preserves the top degree form $\nu^{-2}\in \Lambda^6V^*$, and so lies in $SL(V)$. Now
$$\dim G_0 \le \dim SL(V)=35.$$
But $\dim Spin(12,\C)=66$ so the dimension of the orbit is at least $66-35=31$. On the other hand this orbit lies on the invariant hypersurface $q=3$ which is $32-1=31$-dimensional. The orbit is thus an open subset of the hypersurface and the stabilizer is equal to $SL(V)$.
\vskip .25cm
As we saw,  $\mu(\rho_0)$ is the element $-I/4\in {\mathfrak{gl}}(6,\C)$ which thus acts in the spin representation 
$$S=\Lambda^{ev}V^*\otimes (\Lambda^6 V)^{1/2}$$ as the scalar $(3-p)/4$ in degree $p$. Thus, since $\rho_0=\nu+\nu^{-1}$, the degree $0$ and $6$ elements $\nu$ and $\nu^{-1}$ respectively satisfy 
$$6\nu=4\sigma(\mu(\rho_0))\rho_0+3\rho_0,\qquad 6\nu^{-1}=-4\sigma(\mu(\rho_0))\rho_0+3\rho_0.$$
By invariance $4\sigma(\mu(\rho))\rho+3\rho$ and  $-4\sigma(\mu(\rho))\rho+3\rho$ 
 define, for any $\rho$ on the open orbit, two pure spinors whose sum is $6\rho$. By analyticity, the algebraic expressions above hold for all $\rho$ on the hypersurface $q=3$, and so if $q(\rho)=3$, then $\rho=\alpha+\beta$ where $\alpha$ and $\beta$ are pure spinors with $\langle\alpha,\beta \rangle\ne 0$.
\vskip .25cm
If $q(\rho)\ne 0$, then we rescale to get $q=3$. The choice between  $\alpha$ and $\beta$ is then determined by the choice of square root of $q(\rho)$ so $\alpha$ and $\beta$ are unique up to ordering, which completes the proof.
\vskip .25cm
Now suppose that $\rho$ is real. Proposition \ref{decomprop} says that there are complex pure spinors $\alpha,\beta$ with $\rho=\alpha+\beta$. Reality offers two possibilities: $\alpha$ and $\beta$ are both real, or $\beta=\bar\alpha$. If $\alpha,\beta$ are real then so is $\langle\alpha,\beta\rangle$ and so from (\ref{qform}) $q(\rho)>0$.  If $\beta=\bar\alpha$ then $\langle\alpha,\bar\alpha\rangle$ is imaginary and $q(\rho)<0$. From Proposition \ref{decomprop}, we deduce:
\begin{prp} Let $\rho\in S$ be a real spinor with $q(\rho)<0$. Then $\rho$ is the real part of a pure spinor $\varphi$ with $\langle\varphi,\bar\varphi\rangle\ne 0$.
\end{prp}
These are precisely the pure spinors we need in the definition of a generalized Calabi-Yau manifold. 
\vskip .25cm
When  the vector space $V$ is real, the open set 
$$U=\{\rho\in S: q(\rho)< 0\}$$
is acted on transitively by the real group $\R^*\times Spin(6,6)$. We shall study next the geometry of this space, following closely the parallel discussion of three forms in six dimensions, as in \cite{H1}. In fact, what we are doing here is a direct generalization of that work.

\subsection{The symplectic geometry of the spin representation}\label{sympsection}

\begin{definition} On the open set $U\subset S$ for which $q(\rho)<0$, define the function $\phi$, homogeneous of degree $2$, by
 $$\phi(\rho)=\sqrt{-q(\rho)/3}.$$
\end{definition}
Note from (\ref{qform}) that when we write $\rho=\varphi+\bar\varphi$ for a pure spinor $\varphi$, 
$$i\phi(\rho)=\langle\varphi,\bar\varphi\rangle.$$
\begin{prp} Let $X$ be the Hamiltonian vector field on $U$ defined by the function $\phi$ using the constant symplectic form on $U\subset S$. Describe the vector field on the open set $U$ in the vector space $S$ as a function $X:U\rightarrow S$. Then
\begin{itemize}
\item
 $X(\rho)=\hat\rho$ where $\rho+i\hat\rho=2\varphi$
\item
$X$ generates the circle action $\varphi\mapsto e^{-i\theta}\varphi$
\item
the derivative $DX:U\rightarrow \End S$ defines an integrable almost complex structure $J$ on $U$
\end{itemize}
\end{prp}
\noindent{\bf Proof:} Since $i\phi(\rho)=\langle\varphi,\bar\varphi\rangle$, differentiating along a curve in $U$,
$$i\dot\phi=\langle\dot\varphi,\bar\varphi\rangle+\langle\varphi,\dot{\bar\varphi}\rangle.$$
Up to a scalar, the pure spinors form an orbit, so at each point 
$$\dot\varphi=c\varphi+\sigma(a)\varphi$$
for some $c\in \C$ and $a\in \mathfrak{so}(12,\C)$. But then
\begin{equation}
\langle\dot\varphi,\varphi\rangle=c\langle\varphi,\varphi\rangle+\langle\sigma(a)\varphi,\varphi\rangle=0
\label{doteq}
\end{equation}
where the first term is zero because the bilinear form is skew and the second because, as we saw above, the moment map vanishes on the pure spinors. Using (\ref{doteq}) 
$$\langle\varphi-\bar\varphi,\dot\varphi+\dot{\bar\varphi}\rangle=\langle\dot\varphi,\bar\varphi\rangle+\langle\varphi,\dot{\bar\varphi}\rangle=i\dot\phi.$$
But this can be written as 
$$\dot\phi=\langle\hat\rho,\dot\rho\rangle$$
which means that the Hamiltonian vector field of $\phi$ is $X(\rho)=\hat\rho$.
\vskip .25cm
The circle action in real terms is 
$$\rho\mapsto \cos\theta\rho+\sin\theta\hat\rho$$
so the derivative at $\theta=0$ is $\hat\rho$, the vector field $X$.
\vskip .25cm
Since $\rho+i\hat\rho=2\varphi$, $\hat\rho-i\rho=-2i\varphi$ and so 
$$\hat{\hat\rho}=-\rho.$$
Thus, as a diffeomorphism of $U$, $X\circ X=-id$ and the derivative $J=DX$ thus satisfies $J^2=DX \circ DX=-I$ and defines an almost complex structure on $U$. The proof that it is integrable is the same as in \cite{H1} or \cite{H4} and holds generally for special (pseudo)-K\"ahler manifolds, of which $U$ is an example.

\subsection{The complex structure $J$}\label{Jsection}

The complex structure $J$ on $U\subset S$ turns out to be important in the subsequent development. Recall that $U$ is a homogeneous space of $Spin(6,6)\times \R^*$ under the spin representation. This is a linear action, so every tangent vector to the open set $U$ at $\rho$ is of the form $\sigma(a)\rho$ for some $a$ in the Lie algebra. We show
\begin{prp} \label{J} On the tangent vector $\sigma(a)\rho$, the complex structure $J$ is defined by
$$J(\sigma(a)\rho)=\sigma(a)\hat\rho.$$
Thus the $(0,1)$ vectors are of the form $\sigma(a)\varphi$ where $\rho=\varphi+\bar\varphi$.
\end{prp}
\noindent{\bf Proof:} As $\rho$ varies $\sigma(a)\rho$ defines a vector field $Y$ on $U$. If $a$ is in the Lie algebra of $Spin(6,6)$, then since $\phi$ is invariant and $X$ is the Hamiltonian vector field of $\phi$, we have 
$[X,Y]=0$. 
The central factor $\R^*$ in the group acts by rescaling, so if $a\in \R$  the vector field $Y$ is the Euler vector field -- the position vector $\rho$. Now  $\phi$ is homogeneous of degree $2$ but so is the symplectic form, and this means that $[X,Y]=0$ also.
\vskip .25cm
 Since $J=DX$ and $[X,Y]=0$,
$$J(Y)=DX(Y)=DY(X)=\sigma(a)X=\sigma(a)\hat\rho$$
which proves the proposition.
\vskip .25cm
Although  $J$ is defined on the vector space $S$, it defines a complex structure on the tensor product of $S$ with any vector space and in particular $\Lambda^{ev/od}V^*$, which is where we shall make use of it.
\vskip .25cm
\noindent {\bf Examples:}

\noindent 1. Take the Calabi-Yau case where $\varphi=\Omega$ is a $(3,0)$ form. The space of $(0,1)$-vectors in $\Lambda^{od}V^*\otimes \C$ is from Proposition \ref{J} the image of $\Omega$ under the action of the Lie algebra $\mathfrak{so}(12,\C)+\C$, and using the decomposition ${\mathfrak{so}}(V\oplus V^*)=\End V\oplus \Lambda^2V^*\oplus \Lambda^2V$, this is the $16$-dimensional space of $\Lambda^{od}V^*\otimes  \C$ given by 
$$\Lambda^{3,0}\oplus \Lambda^{2,1}\oplus \Lambda^{3,2}\oplus \Lambda^{1,0}.$$

\noindent 2. In the symplectic case $\varphi=\exp i\omega$, and we obtain for the $(0,1)$ vectors the $16$-dimensional space of  $\Lambda^{ev}V^*\otimes \C$ given by 
$$\exp i\omega\,\C\oplus \exp i\omega(\Lambda^2\otimes\C).$$

\section{The variational problem}

\subsection{The volume functional}
We defined above the function $\phi$ on $U\subset 
\Lambda^{ev/od}V^*\otimes (\Lambda^6V)^{1/2}$. Untwisting by the one-dimensional vector space $(\Lambda^6V)^{1/2}$, there is a corresponding open set, which we still call $U$, in $\Lambda^{ev/od}V^*$ and, since $\phi$ is homogeneous of degree $2$, an invariant function 
$$\phi:U\rightarrow \Lambda^6V^*.$$
The bilinear symplectic form now takes values in  $\Lambda^6V^*$ also and so  the  derivative at $\rho$ of $\phi$ is a linear map from $\Lambda^{ev/od}V^*$ to $\Lambda^6V^*$ which can be written  
\begin{equation}
D\phi(\dot\rho)=\langle\hat\rho,\dot\rho\rangle.
\label{deriv}
\end{equation}
Suppose $M$ is a compact oriented $6$-manifold, and $\rho$ is a form, either odd or even, but in general of mixed degree, which lies at each point of $M$ in the open subset $U$ described above. Following \cite{H2} we shall call such a form {\it stable}. We can then define a volume functional
$$V(\rho)=\int_M\phi(\rho).$$
\begin{thm}\label{var}
A closed stable form $\rho\in \Omega^{ev/od}(M)$ is a critical point of
$V(\rho)$ in its cohomology class if and only if $\rho+i\hat\rho$ defines a \emph{generalized Calabi-Yau structure} on $M$.
\end{thm}
\noindent{\bf Proof:} 
Take the first variation of $V(\rho)$: $$\delta V(\dot\rho)=\int_M
D\phi(\dot\rho)=\int_M\langle\hat\rho,\dot\rho\rangle$$
 from (\ref{deriv}).
The variation is within a fixed cohomology class so
$\dot\rho=d\alpha$. Thus 
$$\delta V(\dot\rho)=\int_M\langle\hat\rho,
d\alpha\rangle=\int_M\sigma(\hat\rho)\wedge d\alpha$$ 
from (\ref{bilinear}).
By Stokes' theorem this is 
$$ \pm\int_M d\sigma(\hat\rho)\wedge \alpha=\pm\int_M \sigma(d\hat\rho)\wedge \alpha=\pm\int_M\langle d\hat\rho,\alpha\rangle$$
since from its definition $\sigma$ commutes with $d$. 

Thus the  variation
vanishes for all $d\alpha$ if and only if $$d\hat\rho=0.$$
A critical point therefore implies $d\varphi=0$ where $2\varphi=\rho+i\hat\rho$.  From Definition \ref{GCYdef} we have a generalized Calabi-Yau manifold.

\subsection{The Hessian}

We shall investigate the Hessian of the functional $V$ at a critical point next. Since $X$ is the Hamiltonian vector field for $\phi$, and $J=DX$ it is clear that $J$ is essentially the second derivative $D^2\phi$. More precisely, we have
\begin{equation}
D^2\phi(\dot\rho_1,\dot\rho_2)=\langle DX \dot\rho_1,\dot\rho_2\rangle=\langle J\dot\rho_1,\dot\rho_2\rangle
\label{symmetric}
\end{equation}
Thus, at a critical point of $V$, the Hessian $H$ is
\begin{equation}
H(\dot\rho_1,\dot\rho_2)=\int_M D^2\phi(\dot\rho_1,\dot\rho_2)=\int_M \langle J\dot\rho_1,\dot\rho_2\rangle
\label{hess}
\end{equation}
where we are restricting the variation to take place in a fixed cohomology class, so that $\dot\rho_1,\dot\rho_2$ are exact forms.
\vskip .25cm
Because of the invariance properties of the functional $V$, any critical point lies on an orbit of critical points, so the Hessian is never non-degenerate. What is the natural group of invariants? 

Firstly $V$ is invariant under diffeomorphisms and those which are homotopic to the identity preserve the de Rham cohomology class of $\rho$ and so the class of forms for the variational problem. The integrand is also invariant under the full group $Spin(6,6)$, so exponentiating sections of the components of the Lie algebra isomorphic to $\Lambda^2T^*$ and $\Lambda^2T$  give further invariant actions. Our variational problem is based on $\rho$ being closed however, and this condition will not be preserved under the action of sections of $\Lambda^2T$. The action of $B\in C^{\infty}(\Lambda^2T^*)$ is the B-field action
$$\rho\mapsto \exp B\wedge\rho.$$
When $B$ is closed, this takes closed forms to closed forms, but to fix the cohomology class we need  $B$ in general to be exact,  for then $B=d\xi$ and 
$$(\exp d\xi)\wedge\rho=\rho+d(\xi\wedge\rho+\frac{1}{2}\xi\wedge d\xi\wedge\rho+\dots)$$ lies in the same cohomology class.

The natural symmetry group of the problem is then the group extension $\mathcal{G}$ 
$$\Omega^{2}_{exact}\rightarrow {\mathcal G}\rightarrow \Diff_0(M).$$
We want to determine when a \GCY   is defined according to Theorem \ref{var} by a Morse-Bott critical point -- non-degenerate transverse to the orbits of the group ${\mathcal G}$. We consider the tangent space to this orbit next.
\vskip .25cm
The action of a vector field on $\rho$ is just the Lie derivative
$${\mathcal L}_X\rho=d\i(X)\rho+\i(X)d\rho=d(\i(X)\rho)$$
since $\rho$ is closed. The infinitesimal action of an exact B-field $B=d\xi$ is
$$d\xi\wedge \rho=d(\xi\wedge\rho).$$
Thus the tangents to an orbit of $\mathcal{G}$  at $\rho$ are  forms
\begin{equation}
\dot\rho=d(\i(X)\rho+\xi\wedge\rho)=d((X+\xi)\cdot\rho)
\label{infaction}
\end{equation}
Because of the invariance of the functional,  if $\alpha$ is exact and $\beta=d(\i(X)\rho+\xi\wedge\rho)$, then $H(\alpha,\beta)=0$. Suppose conversely that the exact form $\beta=d\tau$ has the property that $H(\alpha,\beta)=0$ for all exact forms $\alpha=d\psi$, then from (\ref{hess}),
$$\int_M \langle J d\psi,d\tau\rangle=\pm \int_M \langle \psi,dJd\tau\rangle=0$$
for all $\psi$ so that
$$dJd\tau=0.$$
Thus transverse nondegeneracy is equivalent to the following property:
\begin{definition} \label{ddJ}A  \GCY  is said to satisfy the $dd^J$-lemma if 
$$dJd\tau=0 \Rightarrow d\tau=d(\i(X)\rho+\xi\wedge\rho)$$
for a  vector field $X$ and $1$-form $\xi$.
\end{definition}
This condition may not always be satisfied. Here are two cases when it is:
\begin{prp} The  $dd^J$-lemma holds   if:

\noindent a) the \GCY   is a complex $3$-manifold with a nonvanishing holomorphic $3$-form and which satisfies the $\partial\bar\partial$-lemma, or

\noindent b) it is a symplectic $6$-manifold satisfying the strong Lefschetz condition.

\end{prp}
\noindent{\bf Proof:} 

\noindent 1. When $\rho$ is the real part of a holomorphic $3$-form $\Omega$, $\i(X)\rho=\alpha+\bar\alpha$ where $\alpha$ is of type $(2,0)$ and $\xi\wedge\rho=\beta+\bar\beta$ where $\beta$ is of type $(3,1)$. Moreover $X$ and $\xi$ are uniquely determined by $\alpha$ and $\beta$. Suppose $d\tau$ is an exact odd  form such that $dJd\tau=0$, and write
$$\tau=\gamma_0+\gamma_2+\gamma_4$$ where $\gamma_p$ is a $p$-form. Since we know that $dJd\tau=0$ for $d\tau=d(\i(X)\rho+\xi\wedge\rho)$, we may assume that $\tau$ has no components of type $(2,0)+(0,2)$ and $(3,1)+(1,3)$, so that $\gamma_2\in\Omega^{1,1},\gamma_4\in\Omega^{2,2}$.

From the example above, we saw that $J$ is $-i$ on $\Lambda^{3,0}\oplus \Lambda^{2,1}\oplus \Lambda^{3,2}\oplus \Lambda^{1,0}$, so 
$$Jd\tau=-i(\partial-\bar\partial)\tau$$
and $dJd\tau=0$ implies that for each $p$, $\partial\bar\partial\gamma_p=0$. Now apply the $\partial\bar\partial$-lemma to $d\gamma_p$,  and we see that there  are real forms $\psi_{p-1}$ of degree $p-1$ such that
$$d\gamma_p=i\partial\bar\partial\psi_{p-1}.$$
Here $\psi_{-1}=0$, $\psi_1=\theta_1+\bar\theta_1$ for a $(1,0)$-form $\theta_1$  and $\psi_3=\theta_3+\bar\theta_3$ for a $(2,1)$-form $\theta_3$. 
Thus 
$$d\gamma_2=di(\bar\partial\bar\theta_1-\partial\theta_1)\qquad d\gamma_4=di(\bar\partial\bar\theta_3-\partial\theta_3)$$
and since $\partial\theta_1$ is of type $(2,0)$ and $\partial\theta_3$ is of type $(3,1)$
there exists a real vector field $X$ and a real $1$-form $\xi$ such that 
$$d\tau=d(\gamma_2+\gamma_4)=d(\i(X)\rho+\xi\wedge\rho)$$
as required.
\vskip .25cm
\noindent 2. In the symplectic case, $\rho$ is the real part  $1-\omega^2/2$ of $\exp i\omega$. Thus 
$$d(\i(X)\rho+\xi\wedge\rho)=d(-\i(X)\omega\wedge\omega+\xi-\xi\wedge\omega^2/2).$$
Given an exact even form $d\tau$ with $dJd\tau=0$, write
$$\tau=\gamma_1+\gamma_3+\gamma_5.$$ 
By subtracting $d(\xi\wedge\rho)$ with $\xi=\gamma_1$ we may assume that $\gamma_1=0$. From the example above, $J$ is $-i$ on the forms $\exp i \omega, \exp i\omega\wedge \beta$. To calculate $Jd\tau$ we note that the product with $\omega$ defines an isomorphism from $\Omega^2$ to $\Omega^4$ so the $4$-form $d\gamma_3$ can be written as 
\begin{equation}
d\gamma_3=\omega\wedge\psi
\label{psidef}
\end{equation}
and then the $(0,1)$-part under the action of $J$ is
$$(d\gamma_3)^{(0,1)}=-i\exp i\omega \,\psi$$
and similarly
$$(d\gamma_5)^{(0,1)}=s(i+\omega)\exp i \omega$$
where $d\gamma_5=-2s\omega^3/3$. The equation $dJd\tau=0$ implies $d(d\tau)^{(0,1)}=0$ and this gives
$$d(-i\exp i\omega \,\psi+s(i+\omega)\exp i \omega)=0.$$
The degree $1$ part of this gives $ds=0$ so $s$ is a constant but integrating $d\gamma_5=-2s\omega^3/3$ shows that $s=0$ and hence $d\gamma_5=0$. The degree $3$ part gives $d\psi=0$. Now the strong Lefschetz theorem says that the map $[a]\mapsto [\omega\wedge a]$ from the de Rham cohomology group $H^2(M,\R)$ to $H^4(M,\R)$ is an isomorphism. Since $\psi$ is closed, (\ref{psidef}) shows, using the strong Lefschetz condition, that $\psi$ is exact, and so
$$d\gamma_3=\psi\wedge\omega=d(\xi\wedge\omega)$$
for some $1$-form $\xi$ as required.
\vskip .25cm
Other generalized Calabi-Yau structures satisfy (\ref{ddJ}) because we have:
\begin{prp}
Suppose $(M,\varphi)$ is a $6$-dimensional \GCY   for which  the $dd^J$-lemma  holds. Then the lemma holds for any  transform of $(M,\varphi)$ by a closed B-field.
\end{prp}

\noindent{\bf Proof:} The function $\phi$ is invariant under $Spin(6,6)$ and so in particular under the action $\rho\mapsto\exp B\,\rho$ of the B-field $B$. It follows that the complex structure $J_B$ for the B-field transform is given by
$$J_B=e^ B J e^{-B}.$$ 
So suppose $dJ_Bd\tau=0$. Then since $B$ is closed, 
$$0=d(e^B J e^{-B} d\tau)=e^B d J d (e^{-B}\tau).$$
The $dd^J$-lemma  for the original structure then implies that 
$$d(e^{ -B} \tau)=d(\i(X)\rho+\xi\wedge\rho)$$
so
$$d \tau=d(e^B (\i(X)\rho+\xi\wedge\rho))=d(\i(X)e^B \rho+(\xi-\i(X)B)\wedge e^B\rho)$$
which is of the form
$$d(\i(X)e^B\rho+\eta\wedge e^B\rho)$$
as required.
\vskip .25cm
\noindent{\bf Remark:} By  the definition of $J$,  if $\rho(t)\in U$ is a smooth curve with   tangent vector $\dot\rho$ at $t=0$, the curve $\hat\rho(t)$ has tangent vector $J\dot\rho$ at $t=0$. Thus, globally, if $\rho(t)$ is a smooth family of forms on $M$ defining a generalized Calabi-Yau structure, $d\rho=d\hat\rho=0$ so  the form $\dot\rho$ satisfies
\begin{equation}
dJ\dot\rho=0.
\label{Jinf}
\end{equation}

\subsection{The moduli space}

The non-degeneracy of the critical points implied by the $dd^J$-lemma suggests that the cohomology class of $\rho$ determines the generalized Calabi-Yau structure up to the action of $\mathcal{G}$. We next give a  proof of the existence of a local moduli space as an open set in $H^{ev/od}(M,\R)$. It follows \cite{H1} closely  but because we are forced to be more general, the proof appears simpler. As in \cite{H1} we shall work initially with Sobolev spaces of forms $L^2_k(\Lambda^{ev/od}T^*)$, choosing $k$ appropriately when required.

Choose a metric on $M$ and let $G$ be the Green's operator for the Laplacian $\Delta$ on forms. Then elliptic regularity says that 
$$G:L^2_k\rightarrow L^2_{k+2}$$
and for any form $\alpha$ we have a Hodge decomposition
\begin{equation}
\alpha=H(\alpha)+d(d^*G\alpha)+d^*(Gd\alpha)
\label{Hodge}
\end{equation}
where $H(\alpha)$ is harmonic. 
This $L^2$-orthogonal decomposition is a direct sum of closed subspaces in the  Sobolev space. The exact forms constitute a closed subspace and the $L^2$-orthogonal complement of the image of $d$ is the kernel of $d^*$. Since $J$ is a smooth automorphism of $\Lambda^{ev/od}T^*$, the same is true of the operator $dJ$. From the $dd^J$-lemma we see that the tangent space to the $\mathcal{G}$-orbit through $\rho$ is the closed subspace $\ker dJ$ in the Sobolev space of exact forms.

Choose a transversal to this orbit by taking the $L^2$-orthogonal complement to the tangent space
$$W=\{\alpha \in L^2_k(\Lambda^{ev/od}): d\alpha=0\,\,{\rm and}\,\, 
\alpha\in \im J^*d^*\}.$$

In order to define the volume functional $V$ on $W$ we need  uniform estimates on $\alpha$. From the Sobolev embedding theorem  in $6$ dimensions we can achieve this with $L^2_k$ for $k>3$. Moreover in this range $L^2_k$ is a Banach algebra and so multiplication is smooth. Thus in a Sobolev neighbourhood of $\rho$, the algebraically defined function $V$ is smooth. 

Assigning to a form in $W$ its de Rham cohomology class defines a projection  
$$p:W\rightarrow H^{ev/od}(M)$$ and the derivative of $V$ along the fibres is the linear map on exact forms $d\beta\in W$ given from (\ref{var}) by 
$$\int_M\langle\hat\alpha,
d\beta\rangle=\int_M(\ast\sigma\hat\alpha,d\beta).$$
The critical points of $V$ are thus the zeros of 
$F:W\rightarrow W_{exact}$ defined by
$$F(\alpha)=P(\ast\sigma\hat\alpha)$$
where $P$ is projection onto the exact forms. Its derivative is 
\begin{equation}
DF(\dot\alpha)=P(\ast\sigma J\dot\alpha)
\label{Fderiv}
\end{equation}

\begin{prp}\label{surj}
The derivative of $F$ is surjective at $\rho$.
\end{prp}
\noindent {\bf Proof:} We first need an algebraic lemma:

\begin{lem} \label{Jstar} Take a positive definite inner product on $V$ and put on $\Lambda^{ev/od}V^*$ the induced inner product. Then
$$J^*=\pm\ast\sigma J\sigma\ast$$
the sign depending on whether the forms are even or odd.
\end{lem}

\noindent{\bf Proof:} 
By definition
$$\ast J^*\alpha\wedge \beta=(J^*\alpha,\beta)\nu=(\alpha,J\beta)\nu=\ast\alpha\wedge J\beta$$
where $\nu$ is the volume form. But using (\ref{bilinear}) this can be written in terms of the bilinear form as:
\begin{equation}
\langle \sigma \ast J^*\alpha,\beta\rangle=\langle \sigma \ast\alpha,  J\beta\rangle
\label{Jequation}
\end{equation}
The definition of $J$ in (\ref{symmetric}) gives
$$\langle J\varphi,\psi\rangle =D^2\phi(\varphi,\psi)=\langle J\psi,\varphi\rangle$$
and so from (\ref{Jequation}) we obtain
$$\langle \sigma \ast J^*\alpha,\beta\rangle=\langle J\sigma \ast\alpha, \beta\rangle.$$
Hence 
$$\ast J^*=\sigma J\sigma\ast.$$ 
Since $\ast^2=\pm$ depending on the parity of the forms, we obtain the result.
\vskip .25cm
Now for the proof of Proposition \ref{surj}, suppose $d\psi$ lies in $W$ so there exists $\tau$ such that
$$J^*d^*\tau=d\psi.$$
From Lemma \ref{Jstar} this means that 
$$\sigma\ast J d(\sigma\ast\tau)=d\psi$$
which from (\ref{Fderiv}) says that
$$DF(d(\sigma\ast\tau))=d\psi$$
thus proving the required surjectivity.
\vskip .25cm
The non-degeneracy of the Hessian shows that $DF$ is injective on exact forms, and so from the open mapping theorem and the Banach space implicit function theorem it follows that $F^{-1}(0)$ is a smooth manifold of  forms $\rho$ defining generalized Calabi-Yau structures and which is locally a diffeomorphism under the projection $p$ to an open set in $H^{ev/od}(M,\R)$. This means that an open set in the cohomology space is a local moduli space for these structures -- the modulus of a structure defined by a closed form $\rho$ is just its de Rham cohomology class.
\vskip .25cm
\noindent {\bf Remarks:}

\noindent 1. The fact that the modulus of a generalized Calabi-Yau structure is determined by a cohomology class is not special to six dimensions. For example, Moser's method applied to symplectic forms $\omega$ (see for example \cite{MS}) shows that the cohomology class $[\omega]\in H^2(M,\R)$ is a local modulus up to diffeomorphism for symplectic manifolds and clearly the cohomology class of the B-field in $H^2(M,\R)$ defines $B$ up to exact B-fields. So in particular, for any manifold which  only admits generalized Calabi-Yau structures which are B-field transforms of a symplectic structure -- those of the form 
$\varphi=c\exp (B+i\omega)$
-- the modulus is determined by 
$$(c,[B],[\omega])\in \C^*\times H^2(M)\times H_+^2(M)$$
where $H^2_+(M)$ is the  cone of cohomology classes which carry a symplectic structure. 

\noindent 2. As we saw above, the even type structures in dimension $2$ are all symplectic and in four dimensions  it is only the torus and K3 surface which have  generalized Calabi-Yau structures which are not. In these last two cases we also know enough about the moduli space of ordinary Calabi-Yau structures to say that the generalized local moduli space is an open set in the quadric cone  in $H^{ev}(M,\C)$ defined by $Q(x_0,x_2,x_4)=2x_0\cup x_4-x_2\cup x_2=0$. This is the algebraic condition which the cohomology class of the pure spinor $\varphi$ must satisfy. 
\vskip .25cm
Since a \GCY  is defined by a complex closed form $\varphi$ it is clear that the moduli space, when it exists,  is always  complex. Our description of the six-dimensional case as an open set in the real vector space $H^{ev/od}(M,\R)$ means that this complex structure is less obvious. We describe this next.

\subsection{The special pseudo-K\"ahler structure}

We suppose $\mathcal{U}\subset H^{ev/od}(M,\R)$ is an open set which is a moduli space as above for generalized Calabi-Yau structures on $M$.
\begin{prp}
\label{SK}
 The open set $\mathcal{U}\subset H^{ev/od}(M,\R)$ has the structure of a special pseudo-K\"ahler manifold (see \cite{F}, \cite{H4} or \cite{AL} for a description of this concept). 
\end{prp}
\noindent {\bf Proof:} To define a special pseudo-K\"ahler structure we need  
\begin{itemize}
\item
a flat torsion-free symplectic connection $(\nabla,\omega)$
\item
a complex structure $J$ compatible with $\omega$
\item
a locally defined vector field $X$ such that $\nabla X=J$.
\end{itemize}
In our case the symplectic form is defined by the  constant symplectic form
$$\omega(a,b)=(\sigma(a)\cup b)[M]$$
on the real cohomology 
and we take $\nabla$ to be the flat derivative $D$.

If $\rho$ defines a generalized Calabi-Yau  structure then so does $\cos{\theta}\, \rho+\sin{\theta}\,\hat \rho$ so  assume that $\mathcal{U}$ is invariant by the circle action
$$[\rho]\mapsto \cos{\theta}\, [\rho]+\sin{\theta}\,[\hat \rho].$$
This  action preserves the symplectic form (and has Hamiltonian function given by 
$\Phi=([\rho]\cup[\hat \rho])[M]$, 
 the critical value of the volume functional $V$ at the critical point $\rho$). 

The vector field $X$ generated by this action is $X=[\hat\rho]$ and the derivative of the map $[\rho]\mapsto [\hat\rho]$ on $\mathcal{U}$ defines an almost complex structure $J$ since $\hat{\hat\rho}=-\rho$, and it is integrable just as in \ref{sympsection}. Thus $DX=\nabla X=J$ as required.
\vskip .25cm
\noindent{\bf Remarks:}

\noindent 1.  From the construction of the moduli space we see that every smooth curve in $\mathcal{U}\subset H^{ev/od}(M,\R)$ can be represented by a path of forms $\rho(t)$ defining a generalized Calabi-Yau structure. From (\ref{Jinf}) every cohomology class in $H^{ev/od}(M,\R)$ can therefore be represented by a form $\alpha$ satisfying
$$d\alpha=0,\qquad dJ\alpha=0.$$
The hermitian form of the special K\"ahler metric can then be written in complex form as 
$$([\alpha],[\alpha])=\int_M i\langle \alpha,\bar\alpha\rangle$$
for a closed representative $\alpha$ for which  $J\alpha=i\alpha$. 

From \ref{Jsection}, we can evaluate the hermitian signature of the metric in the two cases of a Calabi-Yau and a symplectic manifold. For the Calabi-Yau threefold, we saw that the $-i$ eigenspace for $J$ was 
$$\Lambda^{3,0}\oplus \Lambda^{2,1}\oplus \Lambda^{3,2}\oplus \Lambda^{1,0}.$$
With the $\partial\bar\partial$-lemma there are closed representatives of each type, and the form is then
$$\int_M \Im(\alpha^{1,0}\wedge\bar\alpha^{3,2})+i\alpha^{3,0}\wedge\bar  \alpha^{3,0}+i\alpha^{2,1}\wedge\bar  \alpha^{2,1}$$
which is of hermitian type $(m,n)$ with
$$m=h^{1,0}+h^{2,1},\qquad n=h^{1,0}+h^{3,0}=h^{1,0}+1.$$
In the symplectic case, the $-i$ eigenspace of $J$ is 
$$\exp i\omega\,\C\oplus \exp i\omega(\Lambda^2\otimes\C).$$
Writing 
$$\alpha=\alpha_0 \exp i\omega+\alpha_1\wedge\exp i\omega$$
where $\alpha_0\in \C$ and $\alpha_1$ is a closed complex $2$-form, we find the hermitian form to be
$$2\int_M \alpha_0\bar\alpha_0 \omega^3-i({\bar\alpha}_0\alpha_1-\alpha_0{\bar\alpha}_1)\wedge\omega^2+\alpha_1\wedge{\bar\alpha}_1\wedge\omega$$
and this has hermitian signature $(\ell,k+1)$ 
where the real quadratic form on $H^2(M,\R)$ defined by 
$$Q(a,b)=(a \cup b\cup [\omega])[M]$$
has signature $(\ell,k)$. This quadratic form  is nondegenerate by the strong Lefschetz property.
\vskip .25cm
\noindent 2. In the symplectic case, the B-field transforms  of $\varphi=\exp {i\omega}$ give an open set in the moduli space and the special pseudo-K\"ahler metric can be seen very concretely. It is known that any special pseudo-K\"ahler metric of dimension $2m$ is determined by a holomorphic function of $m$ variables \cite{F}. In \cite{H4} it was shown that such a description requires a choice of decomposition of the underlying vector space into a direct sum of 
 Lagrangian subspaces. But the even cohomology has a natural decomposition of this type with
$$L=H^0\oplus H^2,\qquad L^*=H^4\oplus H^6$$
both Lagrangian.
To apply the construction of \cite{H4}, we must identify the space of cohomology classes 
$$c\exp([B]+i[\omega])\in H^{ev}\otimes \C=(L\otimes \C)\oplus (L^*\otimes \C)$$
as the Lagrangian submanifold defined by the graph of the derivative of a holomorphic function $\mathcal{F}$ on $L\otimes \C$. For $(c,a)\in H^0(M,\C)\oplus H^2(M,\C)$ this function is 
$$\mathcal{F}=-\frac{1}{6c}a^3[M].$$
The corresponding function in the odd case of the moduli space of Calabi-Yau complex structures is more complicated. Mirror symmetry says that its mirror partner is $\mathcal{F}$ modified by terms from the quantum cohomology of the symplectic manifold $M$.

When $M$ is a torus there are no such corrections, and our special pseudo-K\"ahler structures on the odd or even moduli spaces are just the structures on the open set $U$ in \ref{sympsection}. There we saw that it was only a choice of orientation for $Spin(6,6)$ which differentiated the two moduli spaces. Thus for a torus, although the \emph{objects} parametrized by the two open sets $U$ are very different -- in our language generalized Calabi-Yau structures of odd or even type -- the differential geometric structure of the moduli space is exactly the same.

\subsection{Types of structure in dimension $6$}

Just as we did in \ref{two} and \ref{four} for two and four dimensions we now look at the possible structures in six dimensions which are parametrized by the above moduli spaces. Our approach is ad hoc -- a more systematic treatment of the different algebraic types in all dimensions can be found in \cite{MG}.

To consider the various types, it is useful to note one feature of the geometry of maximally isotropic subspaces -- that two subspaces of even type or odd type intersect in an even-dimensional space, and an odd and even one intersect in an odd-dimensional space. 
\vskip .25cm
Consider first the odd case, so that $\varphi=\varphi_1+\varphi_3+\varphi_5$ is a closed pure spinor. At each point, the isotropic subspace $E$ it defines consists of the complex vectors $v+\xi$ satisfying
\begin{equation}
\iota(v)\varphi_1=0,\quad
\iota(v)\varphi_3+\xi\wedge\varphi_1=0,\quad
\iota(v)\varphi_5+\xi\wedge\varphi_3=0,\quad
\xi\wedge\varphi_5=0. \label{oddpure}
\end{equation}
Now at each point $x\in M$, $E_x\subset(T_x\oplus T_x^*)\otimes\C$ must intersect $T^*_x\otimes \C$ in a subspace of dimension $1$ or $3$, because any higher  dimension will give non-zero vectors in $E_x\cap\bar E_x$. In the first case there is one non-zero $\xi$ with $v=0$ satisfying (\ref{oddpure}) above and five others $v_i+\xi_i$ with $v_1,\dots, v_5$ linearly independent. We extend to a basis $\{v_1,\dots,v_6\}$. Since $\xi\wedge \varphi_i=0$, there exist $\psi_{i-1}$ such that 
$\varphi_i=\xi\wedge\psi_{i-1}$ and then
\begin{equation}
0=\iota(v_i)\varphi_3+\xi_i\wedge\varphi_1=\iota(v_i)(\xi\wedge\psi_2)+\xi_i\wedge\xi\wedge\psi_0=-\xi\wedge(\iota(v_i)\wedge\psi_2+\xi_i\wedge\psi_0)
\label{xidivide}
\end{equation}
using  the isotropic condition on $E$, which  gives
$$0=(v_i+\xi_i,\xi)=-\langle v_i,\xi\rangle=-\i(v_i)\xi.$$
It follows from (\ref{xidivide}) that
$$\psi_2=\sum_1^5\xi_i\wedge\eta_i\wedge\psi_0+\xi\wedge\alpha$$
where $\{\eta_1,\dots,\eta_6\}$ is the dual basis. Similarly, from
$$0=\iota(v_i)\varphi_5+\xi_i\wedge\varphi_3=\iota(v_i)(\xi\wedge\psi_4)+\xi_i\wedge\xi\wedge\psi_2=-\xi\wedge(\iota(v_i)\wedge\psi_4+\xi_i\wedge\psi_2)$$
we obtain
$$\psi_4=\frac{1}{2}\sum_1^5(\xi_i\wedge\eta_i)^2\wedge \psi_0+\xi\wedge\beta$$
and we deduce that 
\begin{equation}
\varphi=\exp (B+i\omega)\varphi_1
\label{type1}
\end{equation}
 for a complex $2$-form $B+i\omega=\sum\xi_i\wedge\eta_i$, well-defined modulo $\varphi_1=\xi\psi_0$. The non-degeneracy condition $\langle \varphi,\bar\varphi\rangle\ne 0$ tells us that $\omega$ is non-degenerate on the $4$-dimensional real subspace of $T_x$ annihilated by the complex $1$-form $\varphi_1$.
\vskip .25cm
The other algebraic type, where $E_x\cap (T^*_x\otimes \C)$ is $3$-dimensional, gives linearly independent $\xi_1,\xi_1,\xi_3$ satisfying (\ref{oddpure}). This means that $\varphi_1=0$ and $\varphi_3=c\xi_1\wedge\xi_2\wedge\xi_3$, the $(3,0)$-form of a complex structure on $T_x$. Then $\varphi_5=(B+i\omega)\wedge\varphi_3$ for some complex $2$-form $B+i\omega$, well-defined modulo $\xi_1,\xi_2,\xi_3$ which means we can take it to be of type $(0,2)$.
\vskip .25cm
These are pointwise descriptions. On a six-manifold with an odd generalized Calabi-Yau structure, the algebraic type may vary from point to point. On an open set where $\varphi_1$ is non-vanishing, we have a local fibration over $\C$ with a symplectic form and B-field along the fibres. If there is any variation, one expects $\varphi_1$ to vanish on a set of measure zero, but if it does vanish on an open set, $\varphi$ defines there an ordinary Calabi-Yau structure, transformed by a complex $2$-form.  The following result gives conditions which do  allow us to determine the algebraic structure at all points.

\begin{prp} Let $M$ be a compact six-dimensional manifold whose first Betti number $b_1(M)$ vanishes.  Then any generalized Calabi-Yau structure $(M,\varphi)$ of odd type on $M$ which satisfies the $dd^J$-lemma  is of the form 
$$\varphi=\Omega+\beta\wedge\Omega$$
where $M$ is a complex manifold with non-vanishing holomorphic $3$-form $\Omega$ and $\beta\in \Omega^{0,2}$ satisfies $\bar\partial\beta=0$. 

If the complex structure on $M$ satisfies the $\partial\bar\partial$-lemma, then $(M,\varphi)$ is equivalent by an exact B-field  to $(M,\Omega)$.
\end{prp}

\noindent{\bf Proof:} let $\varphi=\varphi_1+\varphi_3+\varphi_5$ be the global closed form defining the structure. The algebraic condition of purity on $\varphi$ is preserved not only by scalar multiplication, but also by the action of the scalars in $GL(6,\R)\subset Spin(6,6)$. Thus for each non-zero real number $\lambda$ the following two forms are also closed and pure:
$$\lambda\varphi_1+\lambda\varphi_3+\lambda\varphi_5\qquad \lambda\varphi_1+\lambda^3\varphi_3+\lambda^5\varphi_5.$$
Differentiating at $\lambda=1$, we see that
$$\varphi_1+\varphi_3+\varphi_5,\qquad \varphi_1+3\varphi_3+5\varphi_5$$
are infinitesimal deformations of the structure and hence so is the linear combination
$$\psi=\varphi_1-\varphi_5.$$
From (\ref{Jinf}) we have 
$$dJ\psi=0.$$
However, since $b_1(M)=0$ the de Rham cohomology groups $H^1(M,\R)$ and $H^5(M,\R)$ both vanish, so $\psi$ is exact. The $dd^J$-lemma tells us that 
\begin{equation}
\Re \psi=\Re (\varphi_1-\varphi_5)=d(\i(X)\rho+\xi\wedge\rho)
\label{Xxi}
\end{equation}
for some  vector field $X$ and $1$-form $\xi$. But then $\Re\varphi_1=df$ satisfies
$df=d(\i(X)df)$ and so $f-\i(X)df$ is a constant $c$. 
But this means that $f$ takes the value $c$ at both its maximum and minimum and so must be a constant. Thus  $\Re\varphi_1=df=0$ and similarly $\Im \varphi_1=0$. The algebraic type of $\varphi$ is thus
$$\varphi=\Omega+\beta\wedge\Omega$$
where we can take $\beta$ to be of type $(0,2)$. Since $\varphi$ is closed, $\Omega$ is closed and defines a nonvanishing $(3,0)$-form. Furthermore,  since $\beta\wedge\Omega$  is closed, we have $\bar\partial\beta=0$.
 \vskip .25cm
  Assuming the $\partial\bar\partial$-lemma we have a Hodge decomposition of cohomology so $b_1(M)=0$ implies $h^{0,1}=0$. But  by Serre duality with trivial canonical bundle $h^{0,2}=h^{0,1}$ so $h^{0,2}=0$ and $\beta=\bar\partial \gamma$ for a $(0,1)$-form $\gamma$. Then
$$\beta\wedge\Omega=\bar\partial\gamma\wedge\Omega=d(\gamma+\bar\gamma)\wedge\Omega$$
and $\varphi$ is equivalent under an exact B-field to the Calabi-Yau structure $\Omega$.
\vskip .25cm
We turn now to the  even case, 
$$\varphi=\varphi_0+\varphi_2+\varphi_4+\varphi_6$$
and 
 the dimension of $E_x\cap (T^*_x\otimes \C)$ is $0$ or $2$. In the first case, $E_x$ is the graph of a linear map from $T_x\otimes \C$ to $T_x^*\otimes \C$ which is isotropic and hence given by the vectors 
$$v+\i(v)\beta$$
with $\beta=B+i\omega\in \Lambda^2T^*_x\otimes \C$. In this case
$$\varphi=c\exp (B+i\omega).$$
Since $T_x^*\otimes \C$ is the maximal isotropic subspace defined by an element $\nu\in \Lambda^6T^*_x\otimes \C$, $E_x\cap (T^*_x\otimes \C)\ne 0$ if $\langle\varphi,\nu\rangle\ne 0$ which holds if and only if $\varphi_0\ne 0$.
On a generalized Calabi-Yau manifold the condition $d\varphi=0$ implies that $\varphi_0$ is a constant and so if this algebraic type occurs at one point it occurs everywhere. 

Since $\varphi_2=c(B+i\omega)$ is closed, we see that $B$ and $\omega$  are closed and what we obtain is the B-field transform of a symplectic manifold.
\vskip .25cm
When $\varphi_0=0$, the maximal isotropic subspace is defined by 
\begin{equation}
\iota(v)\varphi_2=0,\quad
\iota(v)\varphi_4+\xi\wedge\varphi_2=0,\quad
\iota(v)\varphi_6+\xi\wedge\varphi_4=0,\quad
\xi\wedge\varphi_6=0. \label{evenpure}
\end{equation}
Since the intersection of $E_x$ with $T^*\otimes \C$ is two-dimensional we now have linearly independent $\xi_1,\xi_2$ satisfying $\xi_i\wedge\varphi_2=0$ which means that $\varphi_2=c\xi_1\wedge\xi_2$. In a similar manner to the odd case we find
$$\varphi=\exp(B+i\omega)\varphi_2.$$
As we saw, this algebraic type  must hold over all the manifold.

 The closed locally decomposable $2$-form $\varphi_2$ defines a codimension $4$ transversally holomorphic foliation and since
$$0\ne \langle \varphi,\bar\varphi\rangle=-\varphi_2\wedge\bar\varphi_4+\bar\varphi_2\wedge\varphi_4=2i\omega\wedge\varphi_2\wedge\bar\varphi_2$$
$\omega$ defines a symplectic structure along the leaves. We also have the  B-field $B$ along the leaves ($B+i\omega$ need only be closed along the leaves). 
A trivial example is the product of a symplectic structure on $S^2$ and a K3 surface.
\vskip .25cm
 For  $6$-manifolds with $b_1(M)=0$ we see then that, under certain generic assumptions, there are essentially  only three types of generalized Calabi-Yau manifold. These are 
\begin{itemize}
\item
Calabi-Yau manifolds
\item
B-field transforms of a symplectic manifold
\item
manifolds with a transversally holomorphic foliation as above.
\end{itemize}

\section{Twisting with a gerbe}

Our variational problem involved an algebraic functional which is invariant under the action of any B-field, closed or not. If we transform the pure spinor $\varphi$ by $B$, then
$e^B\varphi$
is no longer closed, but is so under the modified exterior derivative
$$e^Bd(e^{-B}\alpha)=d\alpha-dB\wedge\alpha.$$
More generally we can replace the exact form $dB$ by a closed $3$-form $-H$ to give an operator on forms:
$$ d_H\alpha=d\alpha+H\wedge\alpha.$$
This differential no longer shifts degrees by one, but maps 
$\Omega^{ev}$ to $\Omega^{od}$ and vice versa. Moreover, our variational approach only distinguishes between even and odd forms, so this differential can be used to set up a similar problem to that in Section 5, but with the form $\rho$ satisfying $d_H\rho=0$  instead of being closed, and now lying in a fixed $d_H$-cohomology class. We obtain the analogue of Theorem \ref{var}:

\begin{thm}\label{Hvar}
A $d_H$-closed stable form $\rho\in \Omega^{ev/od}(M)$ is a critical point of
$V(\rho)$ in its $d_H$-cohomology class if and only if $d_H(\hat\rho)=0$. 
\end{thm}
\noindent{\bf Proof:}
The only different point in the proof is the application of Stokes' theorem. In the even case we have
$$\int_M\sigma(\hat\rho)\wedge (d\alpha+H\wedge \alpha)=\int_M -d\sigma(\hat\rho)\wedge\alpha+H\sigma(\hat\rho)\wedge \alpha$$
but because $H$ is of degree $3$, $H\wedge \sigma(\hat\rho)=-\sigma (H\wedge\hat\rho)$ and so as $d\sigma=\sigma d$,
$$\int_M\sigma(\hat\rho)\wedge d_H\alpha=-\int_M\sigma(d_H\hat\rho)\wedge \alpha.$$
The odd case is similar -- here $H$ anticommutes with the odd form $\sigma(\hat\rho)$.
\vskip .25cm
Given the theorem, it remains to interpret geometrically the condition $d_H(\varphi)=0$ for a pure spinor $\varphi$.  We shall do this in the case where $H/2\pi$ defines an integral class in de Rham cohomology. This corresponds to the appearance of the differential $d_H$ above in recent work on  twisted $K$-theory over the reals (see \cite{A}). We return  to the Courant bracket,  and for motivation, first to the case $p=0$.

Recall the description of the bracket on sections of $T\oplus 1$ in Section 2 as the Lie bracket for invariant vector fields on $M\times S^1$. Instead of the product, 
 suppose $P$ is now a principal $S^1$-bundle over $M$, then  $TP/S^1$ is a vector bundle over $M$ which fits into an exact sequence (the {\it Atiyah sequence}) \cite{A1}
$$0\rightarrow 1\rightarrow TP/S^1\rightarrow TM\rightarrow 0.$$
  Sections of $TP/S^1$ are $S^1$-invariant vector fields on $P$ and so the usual  Lie bracket is defined on these. A \emph{connection} is a splitting of the above exact sequence
$$TP/S^1\cong T\oplus 1$$
  and so, given a connection, a section can be written as $X+f$. A simple calculation shows that the natural Lie bracket on sections can now be written in terms of the Courant bracket as 
  $$[X+f,Y+g]+\i(X)\i(Y)F$$
  where $F\in \Omega^2$ is the curvature of the connection. Thus a flat connection identifies the Courant bracket for $p=0$ with this natural Lie bracket.

A gauge transformation $g:M\rightarrow S^1$ changes a flat connection  by the closed $1$-form $A$ where $iA=g^{-1}dg$ and this is an example of an automorphism of the bracket, defined by a closed form. Note that if $g:M\rightarrow S^1$ is a gauge transformation, the de Rham cohomology class of $A/2\pi$ is integral, so that ``trivial automorphisms" are slightly more general than exact $1$-forms.
\vskip .25cm
Now recall the definition of a \emph{gerbe} \cite{H3} using an open covering $\{U_\alpha\}$. If a gerbe is defined by the cocycle $g_{\alpha\beta\gamma}:U_\alpha\cap U_\beta\cap U_\gamma\rightarrow S^1$ a connection is defined by $1$-forms $A_{\alpha\beta}$ on $U_{\alpha}\cap U_{\beta}$ and $2$-forms $F_\alpha$ on $U_{\alpha}$ such that on the relevant intersections 
\begin{eqnarray*}
F_\beta-F_\alpha&=&dA_{\alpha\beta}\\
iA_{\alpha\beta}+iA_{\beta\gamma}+iA_{\gamma\alpha}&=&g^{-1}_{\alpha\beta\gamma}dg_{\alpha\beta\gamma}
\end{eqnarray*}
 where $H=dF_{\alpha}=dF_\beta$ is a globally defined $3$-form -- called the curvature of the gerbe connection -- such that $H/2\pi$ has an integral de Rham cohomology class. Applying $d$ to the last relation gives
\begin{equation}
dA_{\alpha\beta}+dA_{\beta\gamma}+dA_{\gamma\alpha}=0.
\label{gerbecocycle}
\end{equation}
Over  each open set $U_\alpha$ take $T\oplus T^*$ and identify on the overlap $U_{\alpha}\cap U_{\beta}$ by
$$X+\xi\mapsto X+\xi+\i(X)dA_{\alpha\beta}.$$
Then the cocycle condition (\ref{gerbecocycle}) tells us that we have constructed a vector bundle $W$ which is an extension
$$0\rightarrow T^*\rightarrow W\rightarrow T\rightarrow 0$$
and moreover, since $dA_{\alpha\beta}$ is closed, the Courant bracket is preserved by the identifications. Thus  sections of $W$ inherit a bracket structure. The equation $F_\beta-F_\alpha=dA_{\alpha\beta}$ defines a splitting of the above sequence and hence an isomorphism $W\cong T\oplus T^*$. The bracket operation can now be written using the Courant bracket as
 $$[X+\xi,Y+\eta]+\i(X)\i(Y)H.$$
 We can therefore twist our definition of a generalized complex structure on $M$ by a gerbe -- asking for a maximally isotropic subbundle of $(T\oplus T^*)\otimes \C$ whose sections are closed  under the modified Courant bracket above. A $d_H$-closed form $\varphi=(\rho+i\hat\rho)/2$ arising from the twisted variational problem above then becomes a \emph{twisted generalized Calabi-Yau manifold}.
\vskip .25cm
One final point in this comparison between $p=0$ and $p=1$ for the Courant bracket is that a closed B-field $B$ such that $B/2\pi$ has integral periods should be considered as a gauge transformation. This would mean that the mapping class group for the moduli space of generalized Calabi-Yau structures should be not just the group $\pi_0(\Diff (M))$ of components of $\Diff(M)$ but its extension by $H^2(M,\Z)$. There is evidence both in the physics literature and the symplectic literature \cite {BR} that this integral structure is present.

\end{document}